\documentclass[a4paper,11pt]{article}

\usepackage[utf8]{inputenc}
\usepackage[british]{babel}
\usepackage[margin=2.5cm]{geometry}
\usepackage[shortlabels]{enumitem}
\usepackage{amsmath, amsthm, amssymb, amsfonts, amsrefs}
\usepackage{mathtools}
\usepackage{thmtools, thm-restate}
\usepackage{hyperref}
\usepackage{stmaryrd}
\usepackage{bbm}
\usepackage{tikz}
\usepackage{verbatim}

\declaretheorem[parent=section, name=Theorem]{thm}

\declaretheorem[sibling=thm]{claim}

\definecolor{RoyalAzure}{rgb}{0.0, 0.22, 0.66}
\definecolor{ForestGreen}{rgb}{0.13, 0.55, 0.13}

\DeclareMathOperator*{\Var}{Var}

\hypersetup{
  colorlinks,
  linkcolor={red!60!black},
  citecolor={green!50!black},
  urlcolor={blue!80!black}
}

\title{A simplified proof of the Johansson-Molloy Theorem using the Rosenfeld counting method}

\author{Anders Martinsson\footnote{ETH Z\"{u}rich, maanders@inf.ethz.ch}}

\begin{document}
\maketitle
\begin{abstract}
We show that any triangle-free graph with maximum degree $\Delta$ has chromatic number at most $\left(1+o(1)\right)\Delta/\log \Delta.$ 
\end{abstract}

\section{Introduction}

The Johansson-Molloy Theorem states that any triangle-free graph $G$ with maximum degree $\Delta$ has chromatic number at most $(1+o(1))\Delta/\ln \Delta$ as $\Delta\rightarrow\infty$. In an interesting new development on this problem, a recent article by Hurley and Pirot \cite{PH} present an elementary proof based on the Rosenfeld counting method \cite{Ro}. In fact, their result is quite a bit stronger in that it works even if the underlying graph contains a moderate amount of triangles, it extends to list colorings, and give explicit error terms. In this note, we will present a simplified version of their proof in the case of triangle-free graphs. We prove the following.

\begin{thm}\label{thm:JM} For any $\varepsilon>0$ there exists a constant $\Delta_0=\Delta_0(\varepsilon)$ such that any triangle-free graph with maximum degree at most $\Delta\geq \Delta_0$ is $\lceil (1+\varepsilon) \Delta /\ln \Delta \rceil$-colorable.
\end{thm}

Let us start by fixing some notation. Let $k=\lceil (1+\varepsilon) \Delta /\ln \Delta \rceil$ denote the number of colors, and for any graph $G$, let $\mathcal{C}(G)$ the set of proper $k$-colorings of $G$.

We will show this theorem following the format of Rosenfelds counting method, by proving by induction that the number of proper colorings of a triangle-free graph grows by a large factor whenever a vertex is added. To this end, let $\ell=\ell(\Delta)$ be any function such that $\ln \Delta \ll \ell = \Delta^{o(1)}$, for instance $\ell(\Delta)=\ln^2\Delta.$

\begin{claim}\label{claim} Let $\Delta$ be sufficiently large and let $G$ be any triangle-free graph with maximum degree at most $\Delta$. Then $|\mathcal{C}(G)|/|\mathcal{C}(G-v)|\geq \ell$ for any $v\in G$.
\end{claim}

Note that this immediately implies Theorem \ref{thm:JM}. In fact, the number of proper $k$-colorings of $G$ grows as $\ell^{|V(G)|}>0$.

We will prove Claim \ref{claim} by induction of the number of vertices in $G$. Note that the statement is trivially true when $G$ consists of a single vertex. So we may assume below that $G$ has at least two vertices and that $v\in G$ denotes an arbitrary but fixed vertex.

For any partial proper coloring $c$ of $G$, let $L_c(u)$ denote the set of \emph{available colors} for $u$. That is, the set of colors not present in the neighborhood of $u$. In other words, $L_c(u)$ is the set of colors we can (re-)assign to $u$ while keeping the coloring proper. Note that this definition holds \emph{regardless} of whether $u$ is already assigned a color.

The core of our argument is to study the behavior of a uniformly chosen proper coloring $\mathbf{c}$ of the graph $G-v$. That is, $\mathbf{c}$ is taken uniformly at random from $\mathcal{C}(G-v)$. We first note that the conclusion of the theorem can be expressed in terms of the expectation of the number $|L_\mathbf{c}(v)|$ of available colors for $v$ in $\mathbf{c}$. Observe that each proper coloring of $G$ can be formed by taking a proper coloring of $G-v$ and extending it by coloring $v$ in one of its available colors. Hence, we have the equality
 $$\mathbb{E}[|L_\mathbf{c}(v)|] = \frac{\sum_{c\in \mathcal{C}(G-v)} |L_c(v)| }{|\mathcal{C}(G-v)|} = \frac{|\mathcal{C}(G)|}{|\mathcal{C}(G-v)|},$$
which means that the claim follows if we can prove that, for sufficiently large $\Delta,$ $\mathbb{E}|L_\mathbf{c}(v)| \geq \ell.$ In fact, we will show the stronger statement that, whp as $\Delta\rightarrow\infty$, we have $|L_\mathbf{c}(v)| \gg \ell$.

The proof of this is based on two observations about the distribution of $\mathbf{c}$. Let $u$ be any neighbor of $v$ in $G$, and let $t\geq 0$ be a number to be chosen later. By the induction hypothesis, we know (in fact, for any vertex) that $|\mathcal{C}(G-v)|\geq \ell \cdot |\mathcal{C}(G-v-u)|.$ We will use this to show that it is unlikely for $u$ to have only few available colors in $\mathbf{c}$. Note that any proper coloring of $G-v$ in which $u$ has at most $t$ available colors can be formed by taking some proper coloring of $G-v-u$ and coloring $u$ by picking one out of at most $t$ options. Hence, there are at most $t \cdot |\mathcal{C}(G-v-u)|$ such colorings. Plugging in the induction hypothesis as above, it follows that 
$$\mathbb{P}(|L_\mathbf{c}(u)| \leq t) \leq \frac{t\cdot |\mathcal{C}(G-v-u)|}{|\mathcal{C}(G-v)|} \leq \frac{t}{\ell}.$$
In particular, this means that the expected number of neighbors of $v$ with at most $t$ available colors in $\mathbf{c}$ is at most $t\Delta/\ell$. So if we let $1 \ll t \ll \ell/\ln \Delta$, it follows by Markov's inequality that, with high probability, all but at most $o(k)=o(\Delta/\ln \Delta)$ neighbors of $v$ have more than $t =\omega(1)$ available colors. Crucially, then only a small fraction of colors will be unavailable for $v$ due to these vertices. Thus if suffices to show that it is unlikely for the remaining neighbors of $v$ to block too many of the other $k-o(k)$ colors.

Second, letting $G_0:=G-v-N(v)$, we consider the distribution of $\mathbf{c}$ conditioned on $\mathbf{c}\vert_{G_0} = c_0$ for some $c_0\in \mathcal{C}(G_0)$. Observe that, as $G$ is triangle-free, there cannot be any edges between two vertices in $N(v)$. Thus, for any $u\in N(v)$, we have that $L_\mathbf{c}(u)$ is completely determined by $c_0$ according to $L_{c_0}(u)$. Moreover, any $c\in \mathcal{C}(G-v)$ such that $c\vert_{G_0}=c_0$ can be constructed from $c_0$ by, for each $u\in N(v)$, color $u$ by an arbitrary color in $L_{c_0}(u)$. As $\mathbf{c}$ is uniformly distributed, any such extension of $c_0$ is equally likely. Hence, conditioned on $\mathbf{c}\vert_{G_0}=c_0$, we get that the colors of the neighbors of $v$ are conditionally independent and uniformly chosen from the respective sets $L_{c_0}(u)$. It remains to show that when $L_\mathbf{c}(u)$ is large for most $u\in N(v)$, the set of colors present in $N(v)$ behaves like a coupon collector process.

We have now reduced the problem to the following elementary statement: We are given sets $L_1, L_2, \dots, L_d \subseteq [k]$ for some $d\leq \Delta$ where all but $o(k)$ sets satisfies $|L_i| > t=\omega(1),$ and random variables $X_1\in L_1$, \dots, $X_d \in L_d$ chosen independently and uniformly from the respective sets. We wish to show that, whp, $\mathbf{X}:=[k]\setminus \{X_1, X_2, \dots X_d\}$ contains $\gg \ell$ elements.

Note first that the events $j\in \mathbf{X}$ for $j\in [k]$ are negatively correlated. This is intuitively clear as if some number $j$ never appears among $X_1, \dots, X_d$, then any $X_i$ such that $j\in L_i$ is consequently more likely to hit any value in $L_i\setminus \{j\}$. It follows that $\Var(|\mathbf{X}|)\leq \mathbb{E}|\mathbf{X}|$, so by Chebyshev's inequality it suffices to show that $\mathbb{E}|\mathbf{X}| \gg \ell$.

Condition on the values of $X_i$ for each $i$ such that $|L_i|\leq t$, and let $B$ be the set of colors assigned to these variables. Then $|B|=o(k).$ By linearity of expectation, the conditional expected value of $|\mathbf{X}|$ given these values is
$$\sum_{j\in [k]\setminus B} \prod_{\substack{L_i \ni j\\ |L_i|>t}} \left(1-\frac{1}{|L_i|} \right).$$
It only remains to find a natural lower bound to this expression. By applying the AM-GM inequality to the right-hand side sum, we can lower bound the previous expression by
$$(k-|B|)\left( \prod_{j\in [k]\setminus B} \prod_{\substack{L_i \ni j\\ |L_i|>t}} \left(1-\frac{1}{|L_i|} \right)\right)^{1/(k-|B|)}.$$
Now, by reordering the products, we get $$\prod_{i:|L_i|>t} \prod_{j \in L_i\setminus B} \left(1-1/|L_i|\right) \geq \left( (1-1/t)^t\right)^\Delta = e^{-\left(1+O(1/t)\right)\cdot \Delta},$$ meaning we can further lower bound the above expression by
$$(k-|B|) \exp\left( -  \frac{\left(1+o(1)\right)\cdot \Delta}{k-|B|}\right)= \Theta\left(\frac{\Delta}{\ln \Delta}\right) \cdot \exp\left(-\frac{1+o(1)}{1+\varepsilon} \ln \Delta \right) \geq \Delta^{\varepsilon/(1+\varepsilon)-o(1)},$$
which is clearly $\gg \ell$, as desired. \qed~\\

We end this note with two short remarks. First, by carefully choosing parameters $\ell$ and $t$ in the argument, it is possible to let $\varepsilon$ depend on $\Delta$ according to $\varepsilon=\Omega(\log\log\Delta/\log \Delta)$. Second, it is strictly speaking not needed to argue for concentration of $|\mathbf{X}|$ in the proof. As one only wishes to conclude that $|L_\mathbf{c}(v)|$ is large in expectation, estimating $\mathbb{E}|\mathbf{X}|$ suffices on its own. Nevertheless, we have chosen to write the proof with concentration as we find the statement that $|L_\mathbf{c}(v)|$ is large whp intuitively simpler than that $\mathbb{E}\Big[|L_\mathbf{c}(v)\,|\Big\vert\, \mathbf{c}\vert_{G_0}\Big]$ is large whp.

\section*{Acknowledgements} The author thanks François Pirot and Eoin Hurley for helpful feedback.

\end{document}